\documentclass[a4paper,11pt]{article}

\usepackage{amssymb,amsmath,amsfonts}
\usepackage{array,multirow,ulem,dsfont}
\usepackage{mathtools,stmaryrd,enumitem}
\usepackage{color}
\usepackage{cases}

\newenvironment{proof}{\textit{Proof.}}{\begin{flushright}$\Box$\end{flushright}}

\newtheorem{theorem}{Theorem}
\newtheorem{definition}[theorem]{Definition}

\newtheorem{lemma}[theorem]{Lemma}

\newtheorem{rem}[theorem]{Remark}
\newtheorem{notation}[theorem]{Notation}

\newcommand{\IN}{\mathds{N}}

\newcommand{\IR}{\mathds{R}}

\newcommand{\up}[1]{\textsuperscript{#1}}

\renewcommand{\geq}{\geqslant}
\renewcommand{\leq}{\leqslant}

\begin{document}



\begin{center}
\Large{\bf Rescaling method for blow-up solutions of
nonlinear wave equations}
\end{center}

\begin{center}
Mondher Benjemaa\up{a}, Aida Jrajria\up{b}, Hatem Zaag\up{c}
\end{center}
{\scriptsize \up{a}Faculty of Sciences of Sfax, Sfax University, Tunisia. mondher.benjemaa@fss.usf.tn}\\
{\scriptsize \up{b}Faculty of Sciences of Sfax, Sfax University, Tunisia. aida.jrajria@gmail.com (Corresponding author)}\\
{\scriptsize \up{c}University Sorbonne Paris Nord, LAGA-CNRS, F-93420, Villetaneuse, France. Hatem.Zaag@math.cnrs.fr}\\

\textbf{Abstract:} We develop a hybrid scheme based on a finite difference scheme and a rescaling technique to approximate the solution of nonlinear wave equation. In order to numerically reproduce the blow-up phenomena, we propose a rule of scaling transformation, which is a variant of what was successfully used in the case of nonlinear parabolic equations. A careful study of the convergence of the proposed scheme is carried out and several numerical examples are performed in illustration.\\





\section{Introduction}\label{sec_intro}
In this paper, we are concerned with the study of the numerical approximation of solutions of the nonlinear wave equation that achieve blow-up in finite time
\begin{align}\label{equation1}
\partial_{tt} u = \partial_{xx} u + F(u),\  x\in\,(0,1),\ t \in\,(0,\infty) 
\end{align}
with $F(u) = u^p, p>1$, subject to periodic boundary conditions
\begin{align}\label{boundary condition}
u(0,t) = u(1,t) ,\ t\geq 0 
\end{align}
and the initial conditions
\begin{align}\label{initial condition}
u(x,0) = u_0(x),\ 
\partial_t u(x,0) = u_1(x),\ x\in\,(0,1)
\end{align}
where $u(t): x\in (0,1)\longmapsto u(x,t)\in \mathbb{R}$ is the unknown function. \\
The existence of solutions of the nonlinear wave equation \eqref{equation1}-\eqref{initial condition} was developed in \cite{Caffarelli, Caffarelli and Friedman},, where the authors gave a full description of the blow-up set. In \cite{Glassey}, Glassey proved that under suitable assumptions on the initial data, the solution u of \eqref{equation1} blows up in a finite time the following sense: there exists $T_\infty < \infty$, called the blow-up time, such that the solution $u$ exists on $[0, T_\infty)$ and
\begin{align*}
|| u(.,t)||_\infty \longrightarrow \infty\  \text{as}\  t\longrightarrow T_\infty
\end{align*}
 Caffarelli and Friedman \cite{Caffarelli}  found that there exists the so-called blow-up curve $t = T(x)$  such that the solution $u(x,t)$ satisfies $|u(x,t)|<\infty$  if and only if $t< T(x)$. The blow-up time is therefore $\inf_x T(x)$. For more theoretical results, the reader can refer e.g. to \cite{Azaiez, Merle and Zaag 2003, Merle and Zaag 2005, Merle and Zaag 2007, Merle and Zaag 2012}.
In the numerical direction, the first work was done by Nakagawa in \cite{Nakagawa}  using an adaptive
time-stepping strategy to compute the blowup finite difference solutions and the numerical blow-up time for the 1D semilinear heat equation (see also \cite{abia, abia and al, Cho2007}). For the numerical approximation of blow-up solutions of hyperbolic
equations, Cho applied Nakagawa's ideas to the nonlinear wave equation \cite{Cho2010}. Later on, his results were generalized in \cite{Azaiez and al, Cho2018, Sasaki and Saito}.  \\
In this paper, we intend to develop the rescaling algorithm proposed first by Berger and Kohn \cite{Berg} in 1988 to parabolic equations which are invariant under a scaling transformation. This scaling property allows us to make a zoom of the solution when it is close to the singularity, still keeping the same equation. The scaling transformation is given by
\begin{align}\label{scal transformation}
u_\lambda(\xi,\tau) = \lambda^{\frac{2}{p-1}} u(\lambda \xi,\lambda \tau), \quad \lambda>0.
\end{align}
Clearly, if $u$ is a solution of \eqref{equation1} then $u_\lambda$ is also a solution of \eqref{equation1}.\\
This paper is written as follows. In the next Section, we present the finite difference scheme and the rescaling algorithm. Section \ref{sec Some properties of the discret scheme} is devoted to the proof of several results in concern with the discrete solution. In Section
\ref{sec Convergence of the scheme}, we prove the main results of this paper namely we establish that the
numerical solution converges toward the exact solution. Finally, we give some illustrative examples in Section \ref{sec Numerical examples}.
\section{The numerical algorithm}
In this section, we derive the rescaling algorithm in combination with a finite difference scheme for the nonlinear wave equation \eqref{equation1}.
\subsection{Finite difference approximation}
We use a second order approximation of both the temporal and the spatial derivative operators. Let $I$ be a positive integer and set  $x_i = i\Delta x$  with $\Delta x = \frac{1}{I}$. For the time discretization, let $\Delta t>0$ be a time step and $n\geq 0$ be a positive integer and set $t^n = n\Delta t$. The finite difference scheme of \eqref{equation1} is defined as follows: for all $n\geq 0$ and $1\leq i\leq I$,
\begin{align}\label{scheme}
&\frac{U_{i}^{n+1} - 2U_i^n + U_i^{n-1}}{\Delta t^2}  =  \frac{U_{i+1}^{n} -2U_{i}^{n} + U_{i-1}^{n}}{\Delta x^2} +  F(U_{i}^{n}),
\end{align}
where $U_{i}^{n}$ denotes the approximation for $u(x_i,t^n)$. We set the CFL number $\text{cfl} = \frac{\Delta t}{\Delta x} = 1$, with the following discrete initial and the periodic boundary conditions 
\begin{equation}\label{ini_boundary_cond}
\left\{
\begin{array}{ll}
U_{i}^{0} & \! \! = u_0(x_i), \\
U_{i}^{1} & \! \! = u_0(x_i) + \Delta t u_1(x_i) + \frac{\Delta t^2}{2\Delta x^2}(u_0(x_{i+1}) -2u_0(x_i)+u_0(x_{i-1})) + \frac{\Delta t^2}{2}(F(u_0(x_i))),\\
U_{0}^{n} & \! \! = U_{I}^{n}, \quad U_{I+1}^{n} = U_{1}^{n}.
\end{array}
\right.
\end{equation}

\begin{notation}We denote $U^n$ for $(U_1^n,\cdots, U_I^n)^T$ and  we set 
\begin{align*}
(U_i^n)_t = \frac{U_i^{n+1} - U_i^n}{\Delta t}, \quad (U_i^n)_{\overline{t}} = \frac{U_i^n - U_i^{n-1}}{\Delta t}
\end{align*}
\begin{align*}
(U_i^n)_{t\overline{t}} = \frac{U_i^{n+1} - 2U_i^n + U_i^{n-1}}{\Delta t^2}, \quad (U_i^n)_{x\overline{x}} = \frac{U_{i+1}^n - 2U_i^n + U_{i-1}^n}{\Delta x^2}.
\end{align*} 
We define the norm $|| U||_\infty  = \max_{1\leq i\leq I} |U_i|$ and we write $U\geq 0$ if  $U_i \geq 0$ for all $1\leq i\leq I$. 
 Let $\lbrace (x_i,t^i,U_i^n)| 1\leq i\leq I, n\geq 0 \rbrace$ be a set of data points,  we associate the function $\mathbf{U}$ which is a piecewise linear approximation in both space and time such that for all $(x,t)\in (x_{i+1},x_i)\times(t^n,t^{n+1})$,
 \begin{align}\label{linear approximation power}
 \mathbf{U}(x,t) &= \frac{1}{\Delta t\Delta x}\Big( U_i^n(x_{i+1}-x)(t^{n+1}-t) + U_{i+1}^n(x-x_i)(t^{n+1}-t)\nonumber\\
 & + U_i^{n+1}(x_{i+1}-x)(t-t^n) + U_{i+1}^{n+1}(x-x_i)(t-t^n)\Big).
\end{align}   
\end{notation}
\subsection{The algorithm}
Now, we study the rescaling method for the system \eqref{scheme}-\eqref{ini_boundary_cond}. The transformation \eqref{scal transformation} is originally due to Berger and Kohn \cite{Berg} and was used successfully for parabolic blow-up problems, see \cite{Ngu17}. To set the rescaling algorithm, let $J\in \IN^*$ and $I=J^2$. We consider the partition $[1,I] = \cup_{j=1}^{J}K_j$, with $K_j = [(j-1)J,jJ]$ and the numerical solution $U_{j}^{n} = {U^n}_{\vert K_j}$. Now, we introduce some notations:
\begin{itemize}
\item[-]$0<\lambda <1$: the scale factor such that $\lambda^{-1}$ is a small positive integer.
\item[-] M: the maximum amplitude before scaling.
\item[-] $u^{(k)}(\xi_k,\tau_k)$ is the kth rescaled solution defined in space time variables $(\xi_k,\tau_k)$. The initial index ($k=0$) corresponds to the real solution ($u^{(k=0)} = u, \xi_0 = x, \tau_0 = t$).
\item[-]$U_i^{n,(k)}$: the approximation of $u^{(k)}(\xi_{k,i},\tau_{k}^n)$.
\end{itemize}
 The numerical solution \eqref{scheme} is updated until the first time step $n_0$ such that $||U_j^{n_0}||_\infty \geq M$ is reached. Then using two time levels and a linear interpolation in time to find out a value $\tau_0^*$ satisfying 
\begin{align}\label{find tauStar}
(n_0-1)\Delta t\leq \tau_0^*\leq n_0\Delta t, \, and\, \,  ||\mathbf{U}_j(.,\tau_0^*)|| = M,
\end{align}
as well as the rescaled interval $(x_{i_0^-},x_{i_0^+})$, with $i_0^-,i_0^+ \in K_j$. More precisely, we find the index $i$ where the solution reaches $M$, then we take  
\begin{align}\label{interval to be rescaled power}
\left\{
\begin{array}{ll}
 i_0^+ = i\ \text{and} \  i_0^- =  i - 1 &\text{if $U_j$ is increasing} \\
 i_0^- = i \ \text{and} \  i_0^+ = i + 1 &\text{if $U_j$ is decreasing} \\
 i_0^+ = i+1\ \text{and} \   i_0^- = i-1&\text{otherwise} \\
\end{array}
\right.
\end{align}
The first rescaled solution $u^{(1)}$ is related to u by 
\begin{align*}
u^{(1)}(\xi_1,\tau_1) = \lambda^{\frac{2}{p-1}}u(\lambda\xi_1,\tau^{*}_0+\lambda\tau_1),
\end{align*}
which is also a solution of equation \eqref{equation1} for $\lambda^{-1}x_{i^-_0}<\xi_1<\lambda^{-1}x_{i^+_0}$ and $0<\tau_1<\frac{T-\tau^*_0}{\lambda}$ with  initial conditions
\begin{align*}
&u^{(1)}(\xi_1,0) = \lambda^{\frac{2}{p-1}}u(\lambda\xi_1,\tau^*_0)\\
&u_{\tau_1}^{(1)}(\xi_1,0 )= \lambda^{\frac{p+1}{p-1}}u_t(\lambda\xi_1,\tau_0^*),
\end{align*}
and the boundary conditions  
\begin{align*}
u^{(1)}(\lambda^{-1}x_{i_0^\pm},\tau_1) = \lambda^{\frac{2}{p-1}}u(x_{i_0^\pm},\tau^*_0 + \lambda\tau_1).
\end{align*}  
The maximum value of $u^{(1)}$ at initial time $\tau_1 = 0$ is 
\begin{align*}
||u^{(1)}(.,0) ||_\infty &= \lambda^{\frac{2}{p-1}}||u(.,\tau^*_0)||_\infty\\
&= \lambda^{\frac{2}{p-1}} M.
\end{align*}
Since $\lambda \in (0,1)$, then, $||u^{(1)}(.,0) ||_\infty < M$, i.e the rescaled solution steps down below the threshold criterion. This is the purpose of the rescaling method. Then, we apply the finite difference method to $u^{(1)}$. Let $I_1^{\pm} = \lambda^{-1}i_0^{\pm}$ and $U^{n,(1)}$ the approximation of $u^{(1)}$ at time $\tau_1^n$. Then, the scheme \eqref{scheme} applied to $U^{n,(1)}$ writes: for all $n\geq 0$ and $I_1^-\leq i\leq I_1^+$
\begin{align}
& (U_{i}^{n,(1)})_{t\overline{t}} = (U_{i}^{n,(1)})_{x\overline{x}} + F(U_{i}^{n,(1)}) \nonumber\\
& U_{I_1^-}^{n,(1)} = \psi^{n,(1)},\nonumber\\
& U_{I_1^+}^{n,(1)} = \Psi^{n,(1)},\nonumber\\
& U_{i}^{0,(1)} = \phi^{(1)}_i,\nonumber\\ 
&U_{i}^{1,(1)} = \Phi^{(1)}_i,\nonumber\\
\end{align} 
where 
\begin{align}\label{conditions for U1 power}
&\psi^{n,(1)} = \lambda^{\frac{2}{p-1}}\mathbf{U}(x_{i_0^-},\tau^*_0 + \lambda n\Delta t),\nonumber\\
& \Psi^{n,(1)} = \lambda^{\frac{2}{p-1}}\mathbf{U}(x_{i_0^+},\tau^*_0 + \lambda n\Delta t),\nonumber\\
&\phi_i^{(1)} = \lambda^{\frac{2}{p-1}}\mathbf{U}(\lambda \xi_{1,i},\tau^{*}_0),\nonumber \\
&\Phi_i^{(1)} = \lambda^{\frac{2}{p-1}} \mathbf{U}(\lambda\xi_{1,i},\tau^{*}_0) + \lambda^{\frac{p+1}{p-1}}(\mathbf{U}(\lambda\xi_{1,i},\tau_0^* + \Delta t) - \mathbf{U}(\lambda\xi_{1,i},\tau_0^*))\nonumber\\ 
&+\lambda^{\frac{2p}{p-1}}\Big( \frac{\Delta t^2}{2\Delta x^2}(\mathbf{U}(\lambda\xi_{1,i+1},\tau_0^*) - 2\mathbf{U}(\lambda\xi_{1,i},\tau_0^*) + \mathbf{U}(\lambda\xi_{1,i-1},\tau_0^*))\nonumber\\
& +\frac{\Delta t^2}{2} F(\mathbf{U}(\lambda\xi_{1,i},\tau_0^*))\Big).\nonumber\\
\end{align}
When $||U^{n_{1},(1)}||_\infty$ reaches the given threshold value $M$, we determine $\tau^*_1$ and two grid points $\xi_{1,i_{1}^+}, \xi_{1,i_{1}^-}$ where $i_1^-$ and $i_1^+ \in \lbrace I_1^-,\cdots, I_1^+\rbrace$ using \eqref{find tauStar} and \eqref{interval to be rescaled power} respectively. In the interval where $U^{(1)}\geq  M$ the solution is rescaled further, yielding $U^{(2)}$, and so forth. The $(k+1)^\text{th}$ rescaled solution $u^{(k+1)}$ is introduced when $\tau_k$ reaches a value $\tau^*_k$ satisfying 
\begin{align*}
(n_k-1)\Delta t\leq \tau^*_k\leq n_k\Delta t, \ n_k>0\ \text{ and} \  
||\mathbf{U}^{(k)}(., \tau^*_k)||_\infty = M.
\end{align*}
The interval $(\xi_{k,i_{k}^-},\xi_{k,i_{k}^+})$ to be rescaled satisfies \eqref{interval to be rescaled power} and the solution $u^{(k+1)}$ is related to $u^{(k)}$ by 
\begin{align}\label{the solution u(k) power}
u^{(k+1)}(\xi_{k+1},\tau_{k+1}) = \lambda^{\frac{2}{p-1}}u^{(k)}(\lambda\xi_{k+1},\tau^*_k + \lambda\tau_{k+1}).
\end{align}
Let $I_{k+1}^\pm = \lambda^{-1}i^\pm_k$, the approximation of $u^{(k+1)}(\xi_{k+1,i},\tau_{k+1}^n)$ denoted by $U_i^{n,(k+1)}$  uses the scheme \eqref{scheme}
with the space step $\Delta x$ and the time step $\Delta t$, which reads
\begin{align}\label{rescaling solution Uk power}
&(U_{i}^{n,(k+1)})_{t\overline{t}} = (U_{i}^{n,(k+1)})_{x\overline{x}} + F(U_{i}^{n,(k+1)}) \nonumber\\
& U_{I_1^-}^{n,(k+1)} = \psi^{n,(k+1)},\nonumber\\
& U_{I_1^+}^{n,(k+1)} = \Psi^{n,(k+1)},\nonumber\\
& U_{i}^{0,(k+1)} = \phi^{(k+1)}_i,\nonumber\\ 
&U_{i}^{1,(k+1)} = \Phi^{(k+1)}_i,\nonumber\\
\end{align}
for all $n\geq 0$ and $i \in\lbrace I_{k+1}^-,\cdots,I_{k+1}^+\rbrace$ where 
\begin{align}\label{conditions for Uk power}
&\psi^{n,(k+1)} = \lambda^{\frac{2}{p-1}}\mathbf{U}^{(k)}(x_{i_k^-},\tau^*_k + \lambda n\Delta t),\nonumber\\
& \Psi^{n,(k+1)} = \lambda^{\frac{2}{p-1}}\mathbf{U}^{(k)}(x_{i_k^+},\tau^*_k + \lambda n\Delta t),\nonumber\\
&\phi_i^{(k+1)} = \lambda^{\frac{2}{p-1}}\mathbf{U}^{(k)}(\lambda \xi_{k+1},\tau^{*}_k),\nonumber \\
&\Phi_i^{(k+1)} =  \lambda^{\frac{2}{p-1}} \mathbf{U}^{(k)}(\lambda\xi_{k+1,i},\tau^{*}_k)+ \lambda^{\frac{p+1}{p-1}}(\mathbf{U}^{(k)}(\lambda\xi_{k+1,i},\tau^{*}_k + \Delta t) - \mathbf{U}^{(k)}(\lambda\xi_{k+1,i},\tau^{*}_k)) \nonumber\\
& \qquad \quad+\lambda^{\frac{2p}{p-1}}\Big( \frac{\Delta t^2}{2\Delta x^2}(\mathbf{U}^{(k)}(\lambda\xi_{k+1,i+1},\tau_k^*) - 2\mathbf{U}^{(k)}(\lambda\xi_{k+1,i},\tau_k^*) + \mathbf{U}^{(k)}(\lambda\xi_{k+1,i-1},\tau_k^*)) \nonumber \\
& \qquad \quad+\frac{\Delta t^2}{2} F(\mathbf{U}^{(k)}(\lambda\xi_{k+1,i},\tau_k^*))\Big).
\end{align}
Previously rescaled solutions are stepped forward independently: $\mathbf{U}^{(k)} $ is stepped forward once every $\lambda^{-1}$ time steps of $\mathbf{U}^{(k+1)}$, $\mathbf{U}^{(k-1)}$ once every $\lambda^{-2}$ time steps of $\mathbf{U}^{(k+1)}$, etc. On the other hand, the values of $\mathbf{U}^{(k)}, \mathbf{U}^{(k-1)}$, etc., must be updated to agree with the calculation of $\mathbf{U}^{(k+1)}$. When a time step is reached such that $\Vert \mathbf{U}^{(k+1)} (., \tau_{k+1}) \Vert_\infty > M$, then it is time for another rescaling. Then, the numerical solution $\mathbf{U}_j(x,t)$ of the rescaling method is defined by: for all $1\leq j\leq J$
\begin{align}\label{numerical solution power}
\mathbf{U}_j(x,t) = 
\left\{
\begin{array}{ll}
& \mathbf{U}^{(0)}(x,t)\\
& \lambda^{\frac{-2}{p-1}}\mathbf{U}^{(1)}(\lambda^{-1}x,\lambda^{-1}(t-t_0))\\
&\vdots \\
& \lambda^{\frac{-2(k-1)}{p-1}}\mathbf{U}^{(k-1)}(\lambda^{-(k-1)}x,\lambda^{-(k-1)}(t-t_{k-2}))\\
& \lambda^{\frac{-2k}{p-1}}\mathbf{U}^{(k)}(\lambda^{-k}x,\lambda^{-k}(t-t_{k-1}))\\
\end{array}
\right.
\end{align}
where $t_k = \sum_{i=0}^{k} \lambda^i \tau_i^*$ and $\mathbf{U}^{(k)}$ is the linear interpolation defined in \eqref{linear approximation power} for $k\geq 1$. 
\begin{definition}
We define the numerical blow-up time of $\mathbf{U}_j$ by 
\begin{align}\label{blow up time power}
T^j  = \lim_{k \longrightarrow \infty} \sum_{l=0}^{k} \lambda^l \tau_l^*.
\end{align}
We say that $U_i^n$ blows up if
$£\lim_{n\longrightarrow \infty}||U^n||_\infty  = \infty$£
\end{definition}
Now, we focus on the convergence of the rescaling method. Let $V^{n} = (V_1^n,V_2^n,\cdots,V_I^n)^T$, then one may write \eqref{rescaling solution Uk power} as: for all $n\geq 0$, $i = 1,\cdots,I$ 
\begin{subequations}\label{discrt non zero dirichlet condition power}
\begin{gather}
\begin{alignat}{1}
& (V_{i}^{n})_{t\overline{t}} = (V_{i}^{n})_{x\overline{x}} + F(V_{i}^{n})\\
& V_{1}^{n} = \psi^n,\\
& V_{I}^{n} = \Psi^n,\\
& V_{i}^{0} = \phi_i,\\
& V_{i}^{1} = \Phi_i.
\end{alignat}
\end{gather}
\end{subequations}
where $\psi^n, \Psi^n, \phi_i\  \text{and}\ \Phi_i$ represented by $\psi^{n,(k)}$, $\Psi^{n,(k)}$, $\phi_i^{(k)}$ and $\Phi_i^{(k)}$ in \eqref{conditions for Uk power}. We can see  from \eqref{numerical solution power} that the numerical solution $\mathbf{U}_j$ is built from  $\mathbf{U}^{(k)}$ which is the solutions of the problem \eqref{rescaling solution Uk power}. Thus, we focus on the study of the following problem with the non-periodic Dirichlet conditions:
\begin{align}\label{contenu non zero condition power}
\left\{
\begin{array}{ll}
 \partial_{tt}v = \partial_{xx} v + F(v),\ & x\in (a,b),\ t>0 \\
 v(a,t) = f(t),\ & t\geq 0\\
 v(b,t) = g(t),\ & t\geq 0\\
 v(x,0) = v_0(x),\ & x\in (a,b)\\
 \partial_t v(x,0) = v_1(x),\ & x\in (a,b),
 \end{array} \right.
\end{align}
where $v(t) : x\in (a,b) \longmapsto v(x,t) \in \IR$.
%
\section{Some properties of the discrete scheme}\label{sec Some properties of the discret scheme}
In this section, we give some lemmas of the discrete scheme \eqref{discrt non zero dirichlet condition power} which will be used later. The first lemma below shows a property of the discrete solution.
\begin{lemma}\label{positivity of V_t}Let $V^n = (V_{1}^{n},V_{2}^{n},\cdots,V_{I}^{n})$ be the solution of \eqref{discrt non zero dirichlet condition power}. Denote $Z_i^n = (V_i^n)_{x\overline{x}}$, suppose that  $Z_i^0\geq 0$ and $Z_i^1\geq \dfrac{1}{2}(Z_{i+1}^0 + Z_{i-1}^0)$. Then we have for all $n\geq 0$ 
\begin{align}\label{Positivity of Vxx}
Z_i^n\geq 0.
\end{align}
\end{lemma}
\begin{proof} We proceed by induction on $n$. Suppose that  \eqref{Positivity of Vxx} is valid for all $i$ and $1\leq k \leq n-1$. Now taking into account that $V_i^n$ is a solution of \eqref{discrt non zero dirichlet condition power}, we have 
\begin{align*}
(Z_i^k)_{x\overline{x}} &= ((V_i^k)_{x\overline{x}})_{x\overline{x}} = (V_i^k)_{t\overline{t}x\overline{x}} - (F(V_i^k))_{x\overline{x}}
\end{align*}
and 
\begin{align*}
(Z_i^k)_{t\overline{t}} &= (V_i^k)_{x\overline{x}t\overline{t}}  = (V_i^k)_{t\overline{t}x\overline{x}}.
\end{align*}
 Therefore, by means of Taylor expansions
 \begin{align*}
 (Z_i^k)_{t\overline{t}} - (Z_i^k)_{x\overline{x}} &= (F(V_i^k))_{x\overline{x}}\\
 & = F'(V_i^k) Z_i^k + F''(\zeta_i^k)\frac{(V_{i+1}^k-V_i^k)^2}{2\Delta x^2} + F''(\xi_i^k)\frac{(V_{i-1}^k-V_i^k)^2}{2\Delta x^2} \geq 0.
\end{align*}
Then 
\begin{align*}
Z_i^{k+1} - Z_{i+1}^k\geq Z_{i-1}^k - Z_i^{k-1} \quad \forall\  i \ \text{and}\ 1\leq k\leq n-1.
\end{align*}
It follows
\begin{align*}
Z_i^n &= \sum_{j=0}^{n-1}(Z_{i+j}^{n-j} - Z_{i+j+1}^{n-1-j}) + Z_{i+n}^0\\
&\geq \sum_{j=0}^{n-1} (Z^1_{i-n+1+2j} - Z^0_{i-n+2j+2}) + Z^0_{i+n}\\
&\geq \frac{1}{2}\sum_{j=0}^{n-1} (Z^0_{i-n+2j} - Z^0_{i-n+2j+2}) + Z^0_{i+n}\\
& = \frac{1}{2}(Z^0_{i+n} + Z^0_{i-n})\geq 0.
\end{align*}
\end{proof}
In the next lemma, we show that the numerical solution of \eqref{discrt non zero dirichlet condition power} is not
bounded.
\begin{lemma}Under the same assumptions of Lemma \ref{positivity of V_t}, the numerical solution $V_i^n$ blows up, i.e. $\lim_{n\longrightarrow\infty} V_i^n = \infty$ for all $i$.
\end{lemma}
\begin{proof} Since $(V_i^n)_{x\overline{x}} \geq 0$, we have by \eqref{discrt non zero dirichlet condition power} 
\begin{align*}
(V_i^n)_{t\overline{t}} \geq F(V_i^n),
\end{align*}
implying 
\begin{align*}
(V_i^n)_t \geq (V_i^{n-1})_t + \Delta t F(V_i^n).
\end{align*}
A induction argument yields
\begin{align*}
(V_i^n)_t^2 &\geq \left((V_i^{n-1})_t + \Delta t F(V_i^n)\right)^2\\
& \geq ((V_i^{n-1})_t)^2  + 2 F(V_i^n) (V_i^n - V_i^{n-1})\\
& \geq ((V_i^0)_t)^2 + 2\sum_{k=1}^{n} F(V_i^k)(V_i^k - V_i^{k-1})\\
& \geq ((V_i^0)_t)^2 + 2\int_{V_i^0}^{V_i^n} z^p dz\\
& = ((V_i^0)_t)^2 + \frac{2}{p+1}((V_i^n)^{p+1} - (V_i^0)^{p+1})\\
& = \frac{2}{p+1} (V_i^n)^{p+1} + K_i ,
\end{align*}
with $K_i =((V_i^0)_t)^2 - \frac{2}{p+1} (V_i^0)^{p+1}$. Thus 
\begin{align*}
(V_i^n)_t \geq \sqrt{\frac{2}{p+1}(V_i^n)^{p+1} + K_i}.
\end{align*}
It follows that  
\begin{align*}
V_i^n &\geq V_i^{n-1} + \Delta t \sqrt{\frac{2}{p+1}(V_i^n)^{p+1} + K_i}\\
& \geq  V_i^{n-1} + \Delta t \sqrt{\frac{2}{p+1}(V_i^0)^{p+1} + K_i}\\
&\geq V_i^{0} +n \Delta t \sqrt{\frac{2}{p+1}(V_i^0)^{p+1} + K_i}.
\end{align*}
This achieves the proof.
\end{proof}
The following lemma is a discrete form of the maximum principle.
\begin{lemma}\label{the maximum principle}
Let $b^n=(b_1^n,b_2^n,\cdots,b_I^n)$ be vector such that $b^n\geq 0$. Let $\Theta^n = (\Theta_i^{n})_{1\leq i \leq I}$ satisfy
\begin{align*}
\frac{\Theta_{i}^{n}-2\Theta_{i}^{n-1} + \Theta_{i}^{n-2}}{\Delta t^2} - \frac{\Theta_{i+1}^{n-1}-2\Theta_{i}^{n-1} + \Theta_{i-1}^{n-1}}{\Delta x^2} - b_i^{n-1}\Theta_{i}^{n-1} \geq 0,\quad & 2\leq i\leq I-1\\
\Theta_{1}^{n}\geq 0,\quad &n\geq 0,\\
\Theta_{I}^{n}\geq 0,\quad &n\geq 0,\\
\Theta_{i}^{0} \geq 0, \quad & 1\leq i\leq I,\\
\Theta_{i}^{1} \geq 0, \quad & 1\leq i\leq I\\
\Theta^1_{i} - \Theta^0_{i+1}\geq 0, \quad & 1\leq i\leq I-1. 
\end{align*} 
Then $\Theta^n \geq 0$ for all $n\geq 0$.
\end{lemma}
\begin{proof} Arguing by contradiction, we assume that there exists $n^* \in \IN$ such that there exists  $i^*$ with $\Theta_{i^*}^{n^*}< 0$, and $\Theta_i^n>0$ for all $0\leq n< n^*$. We have
\begin{align*}
\Theta_{i}^{n} \geq \frac{\Delta t^2}{\Delta x^2} (\Theta_{i+1}^{n-1} + \Theta_{i-1}^{n-1}) - \Theta_i^{n-2} + 2(1-\frac{\Delta t^2}{\Delta x^2})\Theta_i^{n-1} +  \Delta t^2 b_i^{n-1}\Theta_{i}^{n-1}.
\end{align*}
Since $\Delta t = \Delta x$, then 
\begin{align}\label{Inegalite sur Theta}
\Theta_i^n - \Theta_{i+1}^{n-1}\geq \Theta_{i-1}^{n-1} - \Theta_i^{n-2} + \Delta t^2 b_i^{n-1} \Theta_{i}^{n-1}
\end{align}
Let $W_i^n = \Theta_{i}^{n} - \Theta_{i+1}^{n-1}$, it follows from \eqref{Inegalite sur Theta}
\begin{align*}
W_{i}^{n} &\geq  W_{i-1}^{n-1}  + \Delta t^2 b_i^{n-1}\Theta_{i}^{n-1}\\
& \geq W_{i-n+1}^{1} + \Delta t^2 \sum_{l=1}^{n-1} b_{i+1-l}^{n-l}\Theta_{i+1-l}^{n-l}.
\end{align*}
Then 
\begin{align*}
\Theta_{i^*}^{n^*} &= \Theta_{i^*+1}^{n^*-1} + W_{i^*}^{n^*}\\
& = \Theta_{i^*+n^*}^0 + \sum_{k=0}^{n^*-1} W_{i^*+k}^{n^*-k}\\
& \geq \Theta_{i^*+n^*}^0 + \sum_{k=0}^{n^*-1}W_{i^*-n^*+1+2k}^1 + \Delta t^2\sum_{k=0}^{n^*-2}\ \sum_{l=1}^{n^*-k-1} b_{i^*+k+1-l}^{n^*-k-l} \Theta_{i^*+k+1-l}^{n^*-k-l} \\
& \geq \Theta_{i^*+n^*}^0 + \sum_{k=0}^{n^*-1} (\Theta^1_{i^*-n^*+1+2k} - \Theta^0_{i^*-n^*+2+2k}) + \Delta t^2\sum_{k=0}^{n^*-2}\  \sum_{l=1}^{n^*-k-1} b_{i^*+k+1-l}^{n^*-k-l} \Theta_{i^*+k+1-l}^{n^*-k-l}\\
& \geq 0 
\end{align*}
which is a contradiction.
\end{proof}
\section{Convergence of the scheme}\label{sec Convergence of the scheme}
We prove in the following the convergence of the scheme  \eqref{discrt non zero dirichlet condition power}. The next result establishes that for each fixed time interval $[0,T_\infty)$ where the solution \eqref{contenu non zero condition power} $v$ is defined, the numerical solution of the problem \eqref{discrt non zero dirichlet condition power} approximates $v$ as $\Delta x\longrightarrow 0$.
\begin{theorem}\label{theorem of convergence power}
Let $V_{i}^{n}$ and $v$ be the solution of \eqref{discrt non zero dirichlet condition power} and \eqref{contenu non zero condition power} respectively. Let $T_\infty$ denotes the blow-up time of $v$ and let $T_0$ be an arbitrary number such that $0 < T_0 < T_\infty$. Suppose that $v \in C^{2}([a,b]\times [0,T_0])$ and the initial data and boundary data of \eqref{discrt non zero dirichlet condition power} satisfy 
\begin{align*}
\epsilon_1 = \sup_{x\in[a,b]}|v_0(x)- \phi(x)| = o(\Delta x) \quad \text{as} \ \Delta x\longrightarrow 0,
\end{align*}
\begin{align*}
\epsilon_2 = \sup_{x\in[a,b]}|v_1(x)- \Phi(x)| = o(\Delta x) \quad \text{as} \ \Delta x\longrightarrow 0,
\end{align*}
\begin{align*}
\epsilon_3 = \sup_{t\in[0,T_0]}|f(t)- \psi(t)| = o(1) \quad \text{as} \ \Delta t\longrightarrow 0,
\end{align*}
\begin{align*}
\epsilon_4 = \sup_{t\in[0,T_0]}|g(t)- \Psi(t)| = o(1) \quad \text{as} \ \Delta t\longrightarrow 0,
\end{align*}
where $\phi$, $\Phi$, $\psi$ and $\Psi$ are the interpolations of $\phi_i$, $\Phi_i$, $\psi^n$ and $\Psi^n$ respectively defined in \eqref{linear approximation power}. Then, 
\begin{align}\label{convergence of solution power}
\max_{0\leq n\leq N}||V_i^n - v(x_i,t^n)||_\infty = \mathcal{O}(\epsilon_1 + \epsilon_2 + \epsilon_3 + \epsilon_4 + \Delta x^2) \ \text{as}\quad \Delta x\longrightarrow 0,
\end{align}
where $N>0$ such that  $t^N \leq T_0$.
\end{theorem}
\begin{proof}
 We denote $N^*$ the greatest value such that $N^*< N$, and for all $0\leq n < N^*$
 \begin{align}\label{e less 1 power}
 || V_i^n - v(x_i,t^n)||_\infty < 1.
 \end{align}
Let $e^n_i = V_i^n - v(x_i,t^n)$ the error at the node $(x_i,t^n)$. By Taylor's expansion and \eqref{contenu non zero condition power}, we have for all $2\leq i\leq I-1$ and $0\leq n < N^*$
\begin{equation*}
\frac{v(x_i,t^{n+1})-2v(x_i,t^n)+v(x_i,t^{n-1})}{\Delta t^2} = \partial_{tt}v(x_i,t^n) +  \frac{\Delta t^2}{24}\lbrace{\partial_{tttt} v(x_i,\tilde{t^n}) + \partial_{tttt}v(x_i,\tilde{\tilde{t^n}})\rbrace},
\end{equation*}
 where $\tilde{t^n}, \tilde{\tilde{t^n}} \in [t^{n-1},t^{n+1}]$ and 
\begin{equation*}
\frac{v(x_{i+1},t^{n})-2v(x_i,t^n)+v(x_{i-1},t^{n})}{\Delta x^2} = \partial_{xx}v(x_i,t^n) + \frac{\Delta x^2}{24}\lbrace{\partial_{xxxx} v(\tilde{x_i},t^n) + \partial_{xxxx}v(\tilde{\tilde{x_i}},t^n)\rbrace},
\end{equation*} 
 where $\tilde{x_i}, \tilde{\tilde{x_i}} \in [x_{i-1},x_{i+1}]$. Using the mean value theorem, we obtain 
\begin{equation*}
F(V_i^n) - F(v(x_i,t^n)) = F'(\delta_i^n) (V_i^n - v(x_i,t^n)),
\end{equation*}
where $\delta_i^n$ is an intermediate value between $V_i^n$ and $v(x_i,t^n)$. It follows
\begin{equation*}
(e_i^n)_{t\overline{t}} - (e_i^n)_{x\overline{x}} = F'(\delta_i^n) e_i^n + r_i^n, 
\end{equation*} 
with
\begin{align*}
r_i^n = - \frac{\Delta t^2}{24}\lbrace{\partial_{tttt} v(x_i,\tilde{t^n}) + \partial_{tttt}v(x_i,\tilde{\tilde{t^n}})\rbrace} + \frac{\Delta x^2}{24} \lbrace{\partial_{xxxx} v(\tilde{x_i},t^n) + \partial_{xxxx}v(\tilde{\tilde{x_i}},t^n)\rbrace}.
\end{align*}
Let  $C$ be positive constant such that
 \begin{align*}
 \frac{1}{12}\max_{(x,t)\in[a,b]\times[0,T_0]}|\left(\partial_{tttt} v(x,t)| +|\partial_{xxxx} v(x,t)|\right) \leq C.
 \end{align*}
Since $\Delta t = \Delta x$, we obtain for all $2\leq i\leq I-1$ and $n\geq 0$
\begin{align*}
(e_i^n)_{t\overline{t}} - (e_i^n)_{x\overline{x}}  \leq F'(\delta_i^n) e_i^n + C \Delta x^2.
\end{align*} 
Now, consider the function $E(x,t)$ defined by 
\begin{equation*}
E(x,t) = e^{K t + x}(\epsilon_1 + \epsilon_2 + \epsilon_3 + \epsilon_4 + C \Delta x^2) ,
\end{equation*}
with $K$ is a positive constant which will be chosen adequately. Using Taylor expansion, we get 
\begin{align}\label{Taylor on E power}
&(E(x_i,t^n))_{t\overline{t}} - (E(x_i,t^n))_{x\overline{x}} - F'(\delta_i^n) E(x_i,t^n)\nonumber\\
& = \partial_{tt}E(x_i,t^n) - \partial_{xx} E(x_i^n,t) - F'(\delta_i^n) E(x_i,t^n) \nonumber\\
& + \frac{\Delta t^2}{24}\{\partial_{tttt}E(x_i,\bar{t^n}) + \partial_{tttt}E(x_i,\bar{\bar{t^n}}) \} - \frac{\Delta x^2}{24}\{\partial_{xxxx}E(\bar{x_i},t^n) + \partial_{xxxx}E(\bar{\bar{x_i}},t^n)\}
\end{align}
where $\bar{t^n},\bar{\bar{t^n}} \in [t^{n-1},t^{n+1}]$ and $\bar{x_i}, \bar{\bar{x_i}} \in [x_{i-1},x_{i+1}]$. We have for all $x\in [a,b]$ and $t\in [0,T_0]$ 
\begin{equation*}
\partial_{tt}E(x,t) - \partial_{xx} E(x,t) = (K^2-1)E(x,t)
\end{equation*}
and
\begin{equation*}
E(a,0)\leq E(x,t)\leq E(b,T_0),
\end{equation*}
yielding
\begin{align*}
\partial_{tttt} E(x,t)  = K^4 E(x,t) \geq K^4 E(a,0)
\end{align*}
and 
\begin{align*}
\partial_{xxxx} E(x,t)  =  E(x,t) \leq E(b,T_0).
\end{align*}
Then, \eqref{Taylor on E power} implies
\begin{align*}
& (E(x_i,t^n))_{t\overline{t}} - (E(x_i,t^n))_{x\overline{x}}  - F'(\delta_i^n) E(x_i,t^n)\\
 & \geq (K^2 - 1 - F'(\delta_i^n)) E(x_i,t^n) + \frac{\Delta x^2}{12}\{ (K^4 E(a,0) - E(b,T_0)\}.
\end{align*}
By taking $K$ large enough such that the right hand side of the above inequality is large than $C\Delta x^2$, we obtain 
\begin{align*}
(E(x_i,t^n))_{t\overline{t}} - (E(x_i,t^n))_{x\overline{x}}  - F'(\delta_i^n) E(x_i,t^n) - C\Delta x^2 \geq 0.
\end{align*}
Therefore, from Lemma \ref{the maximum principle} with $b_i^n = F'(\delta_i^n)$ and $\Theta_i^n = E(x_i,t^n)- e_i^n$, we get for $1\leq i\leq I$ 
\begin{align*}
E(x_i,0) &\geq e_i^0\\
\partial_t E(x_i,0) &\geq e_i^1\\
E(x_1,t^n) &\geq e_1^n\\
E(x_I,t^n) &\geq e_I^n,
\end{align*}
and
\begin{align*}
E(x_i,\Delta t) - E(x_{i+1},0) &= e^{K \Delta t + x_i}(\epsilon_1 + \epsilon_2 + \epsilon_3 + \epsilon_4 + C \Delta x^2) - e^{x_{i+1}}(\epsilon_1 + \epsilon_2 + \epsilon_3 + \epsilon_4 + C \Delta x^2)\\
 & =  (e^{K \Delta t} e^{x_i} - e^{x_i + \Delta x})(\epsilon_1 + \epsilon_2 + \epsilon_3 + \epsilon_4 + C \Delta x^2)\\
 & =  e^{x_i}(e^{K \Delta t} - e^{ \Delta x})(\epsilon_1 + \epsilon_2 + \epsilon_3 + \epsilon_4 + C \Delta x^2)\\
 & \geq e_i^1,
\end{align*}
yielding
\begin{align*}
E(x_i,\Delta t) - e_i^1 \geq E(x_{i+1},0) - e_{i+1}^0,
\end{align*}
and hence 
\begin{align*}
V_i^n - v(x_i,t^n) \leq E(x_i,t^n). 
\end{align*}
Using the same argument for $\epsilon_i^n =v(x_i,t^n) - V_i^n = - e_i^n$, we obtain
\begin{align*}
(\epsilon_i^n)_{t\overline{t}} - (\epsilon_i^n)_{x\overline{x}} &\leq  F(v(x_i,t^n)) - F(V_i^n) + C\Delta x^2\\
&\leq  F'(\delta_i^n)\epsilon_i^n + C\Delta x^2.
\end{align*}
By Lemma \ref{the maximum principle} with $\Theta_i^n = E(x_i,t^n) - \epsilon_i^n$, we get
\begin{align*}
v(x_i,t^n) - V_i^n \leq E(x_i,t^n),
\end{align*}
then 
\begin{align*}
|V_i^n - v(x_i,t^n)| &\leq E(x_i,t^n)\\
& \leq e^{K T_0 + b}(\epsilon_1 + \epsilon_2 + \epsilon_3 + \epsilon_4 + C \Delta x^2). 
\end{align*}
Hence, we obtain, for $ n < N^*$
\begin{align}\label{inequality on e power}
\max_{1\leq i\leq I}|V_i^n - v(x_i,t^n)| &\leq e^{K T_0 + b}(\epsilon_1 + \epsilon_2 + \epsilon_3 + \epsilon_4 + C \Delta x^2).
\end{align}
In order to prove \eqref{convergence of solution power}, we have to show that $N^* = N$. If it is not true, we have by \eqref{e less 1 power} and \eqref{inequality on e power} 
\begin{align*}
1\leq || V_i^{N^*} - v(x_i,t^{N^*})||_\infty \leq e^{K T_0 + b}(\epsilon_1 + \epsilon_2 + \epsilon_3 + \epsilon_4 + C \Delta x^2).
\end{align*}
The last term of the above inequality goes to zero as $\Delta x$ tends to zero, which is a contradiction. The proof is achieved.
\end{proof}
\begin{rem}
From the relation between the numerical solution and $\mathbf{U}^{(k)}$ in \eqref{numerical solution power}, we conclude the convergence of the rescaling method according to the Theorem \ref{theorem of convergence power}.
\end{rem}
\section{Numerical examples}\label{sec Numerical examples}
In this section, we present some numerical examples. For all the examples, we set $\lambda = \frac{1}{2}$ and we choose the threshold value $M$ such that the maximum of the initial data of all rescaled solutions are equal, i.e. for all $k\geq 0$ we have $|| u_0 ||_\infty = \lambda^{\frac{2}{p-1}}|| u^{(k)}(\tau_k^*)||_\infty$. Since $||u^{(k)}||_\infty = M$, then $M = \lambda^{\frac{-2}{p-1}}|| u_0||_\infty$.
\subsection*{Example 1}\label{example 1}
\begin{figure}
\begin{center}
\includegraphics[width=8cm,height=6cm]{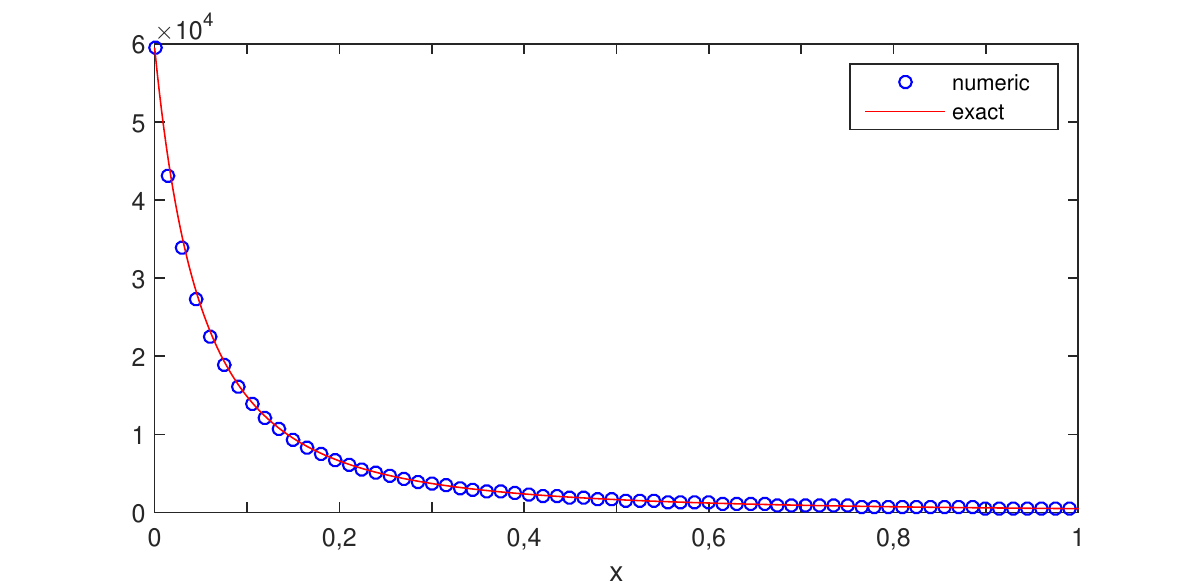}\hspace*{-2mm}
\includegraphics[width=8cm,height=6cm]{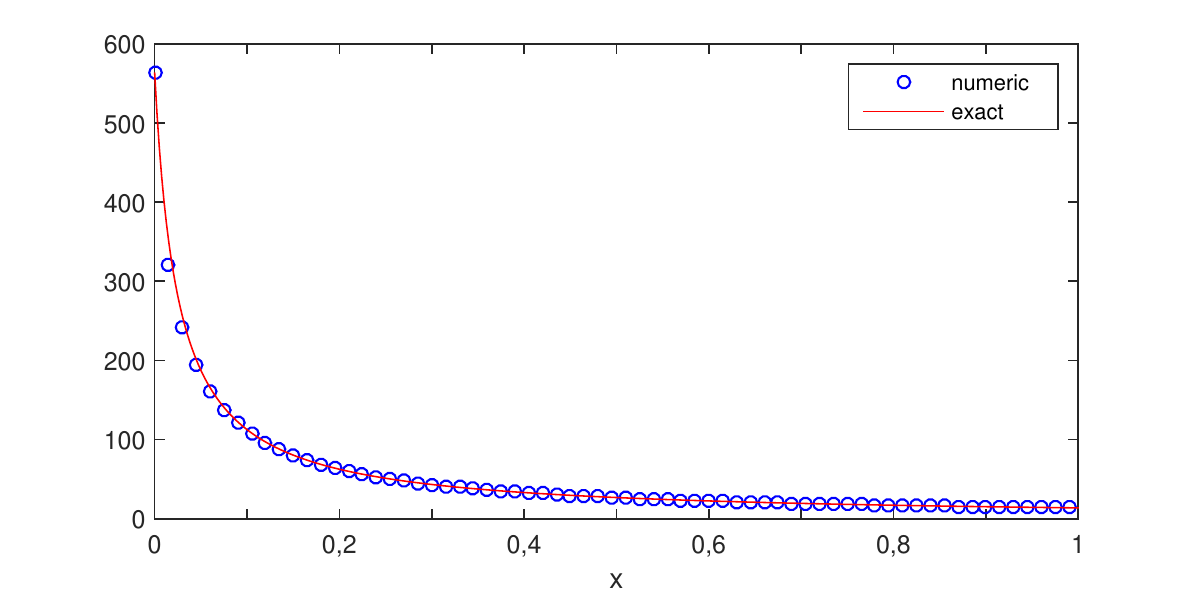}
\caption{\label{comparison_exp1}Exact solution (red line) and numerical solution (blue circles) with p = 2 (left) and p = 3 (right).}
\end{center}
\end{figure}

\begin{figure}
\begin{center}
\includegraphics[width=8cm,height=6cm]{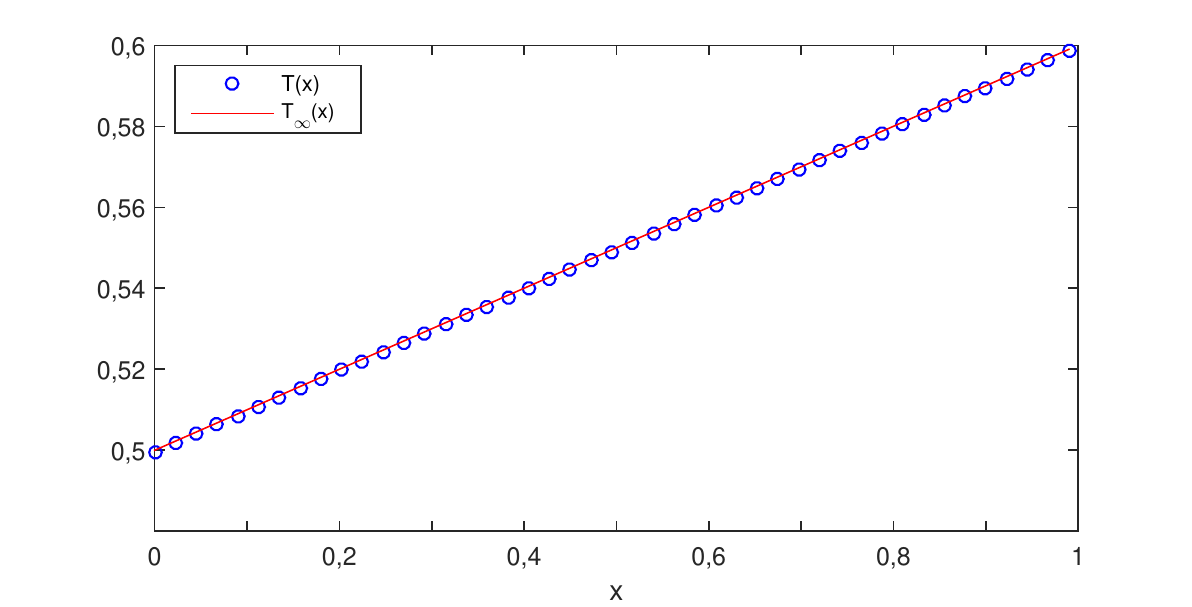}\hspace*{-2mm}
\includegraphics[width=8cm,height=6cm]{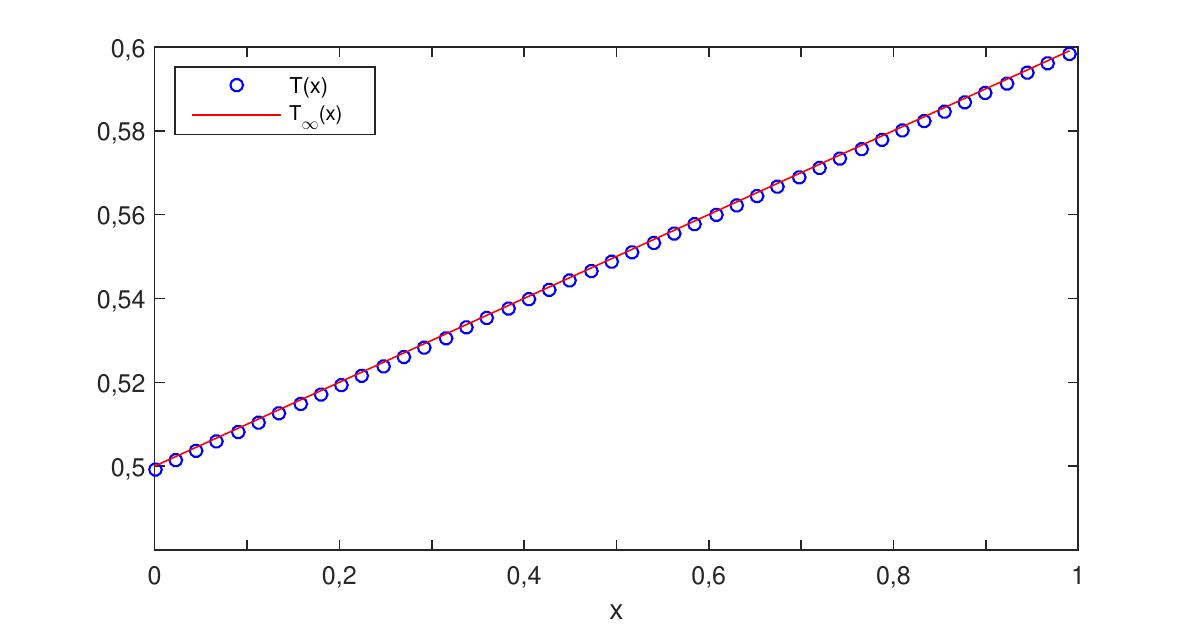}
\caption{\label{blowup_curve_exp1}Comparison between the numerical blowup time (blue circles) and $T_\infty$(red line) for $p = 2$ (left) and $p = 3$ (right).}
\end{center}
\end{figure}

\begin{figure}
\begin{center}
\includegraphics[width=8cm,height=6cm]{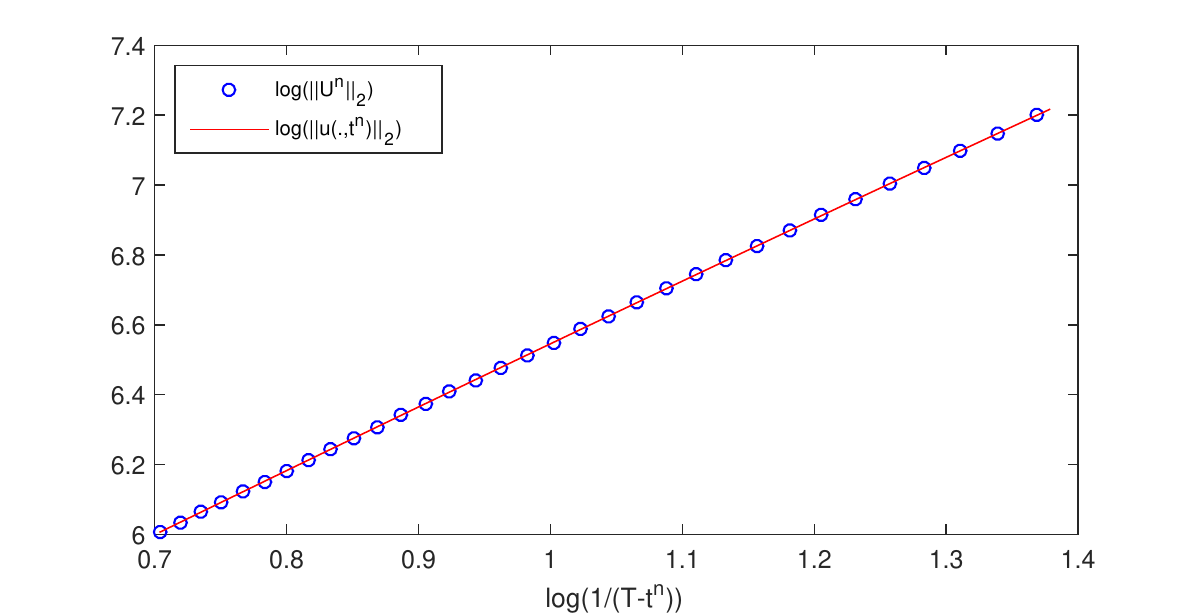}\hspace*{-2mm}
\includegraphics[width=8cm,height=6cm]{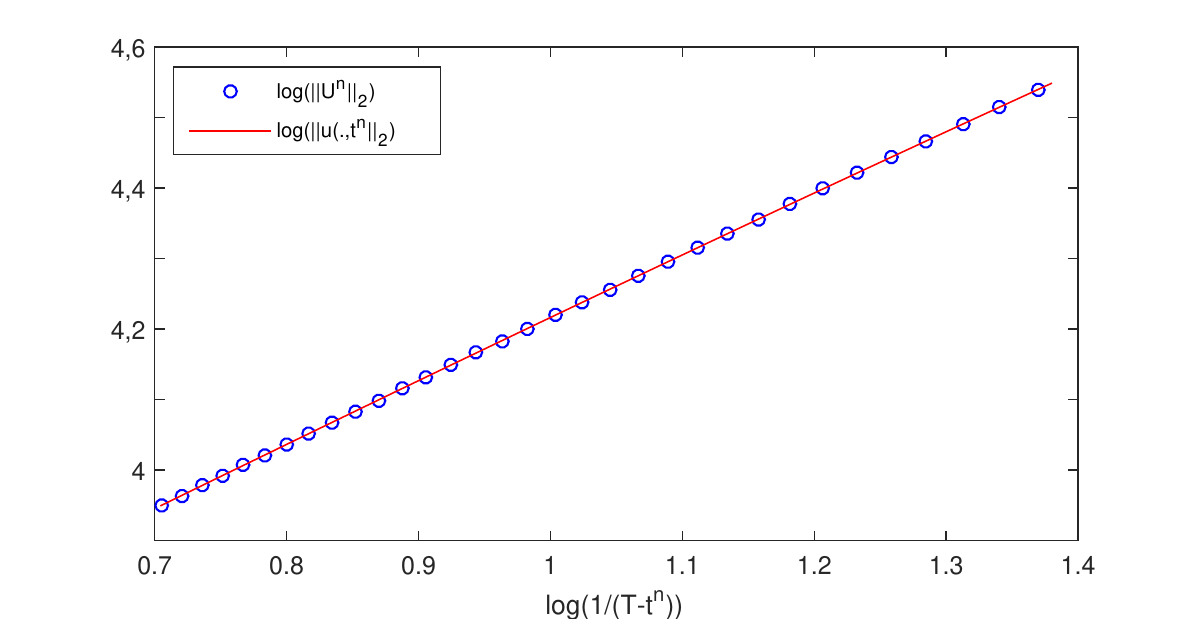}
\caption{\label{blowup_rate_exp1}Blow up rate for $p = 2$ (left) and $p = 3$ (right).}
\end{center}
\end{figure}

We consider the system  \eqref{equation1} with an exact solution given by
\begin{align*}
u(x,t) = \mu(T-t+dx)^{\frac{2}{1-p}}
\end{align*}
with $\mu = \Big(2(1-d^2)\frac{p+1}{(p-1)^2}\Big)^{\frac{1}{p-1}}$ and $d\in (0, 1)$ is an arbitrary parameter. The parameters used are $T = 0,5$ and $d = 0.1$. Figure \ref{comparison_exp1} shows a comparison between the exact solution and the numerical solution. One can notice a very good superposition between the solutions. In table \ref{tab_norm}, we report the relative $L^2$ and $L^\infty$ errors. The blow-up time for both cases is set $T_\infty(x) = T + dx$. Thus, one can approximate numerically the blow-up curve $T_\infty(x)$ by computing the numerical blow up time $T^j$ \eqref{blow up time power} for all $1\leq j\leq J$. Figure \ref{blowup_curve_exp1} shows a comparison between the exact blow up curve and $T(x)$. In \cite{Merle and Zaag 2003}, the authors proved that the solution satisfies
\begin{align*}
|| u(.,t)||_{2} \sim (T-t)^{\frac{-2}{p-1}},
\end{align*}
one has 
\begin{align*}
\log(|| u(.,t)||_{2}) \sim \frac{-2}{p-1}\log(T-t)
\end{align*}
In order to calculate the blow-up rate $\frac{2}{p-1}$ numerically, we show the plot of $\log(|| U^n||_2)$ versus $\log(\frac{1}{T-t^n})$. Figure \ref{blowup_rate_exp1} presents these slopes for $p = 2$ and $p = 3$. 
\renewcommand{\arraystretch}{1.25}
\begin{table}
\begin{center}
\begin{tabular}{|*{3}{c|}}
\hline
\multicolumn{3}{|c|}{$p = 2$}\\ 
\hline 
I & \rule[12pt]{0pt}{5pt} $\frac{\|U_{num}-u_{exact}\|_2}{\|u_{exact}\|_2}$ \rule[-8pt]{0pt}{5pt} & $\frac{\|U_{num} - u_{exact}\|_\infty}{\|u_{exact}\|_\infty}$ \\ 
\hline 
$2^6$  & $9\times 10^{-2}$ & $10\times 10^{-2}$ \\ 
\hline 
$2^7$  & $3\times 10^{-2}$ & $3.2\times 10^{-2}$ \\ 
\hline 
$2^8$  & $8.2\times 10^{-3}$ & $8.7\times 10^{-3}$ \\ 
\hline 
$2^9$  & $2.1\times 10^{-3}$ & $2.2\times 10^{-3}$ \\ 
\hline 
\end{tabular}
\begin{tabular}{|*{3}{c|}}
\hline
\multicolumn{3}{|c|}{$p = 3$}\\ 
\hline 
I & \rule[12pt]{0pt}{5pt} $\frac{\|u_{num}-u_{exact}\|_2}{\|u_{exact}\|_2}$ \rule[-8pt]{0pt}{5pt}& $\frac{\|u_{num}-u_{exact}\|_\infty}{\|u_{exact}\|_\infty}$ \\ 
\hline 
$2^6$  & $7\times 10^{-2}$ & $9\times 10^{-2}$ \\ 
\hline 
$2^7$  & $2.3\times 10^{-2}$ & $3.3\times 10^{-2}$ \\ 
\hline 
$2^8$ & $6.6\times 10^{-3}$ & $9.5\times 10^{-3}$ \\ 
\hline 
$2^9$ & $1.7\times 10^{-3}$ & $2.5\times 10^{-3}$ \\ 
\hline 
\end{tabular}
\caption{\label{tab_norm} Relative errors of the numerical solution  versus the exact solution.}
\end{center}
\end{table}
%
%
\subsection*{Example 2}\label{example 2}
\begin{figure}
\begin{center}
\includegraphics[width=9cm,height=7cm]{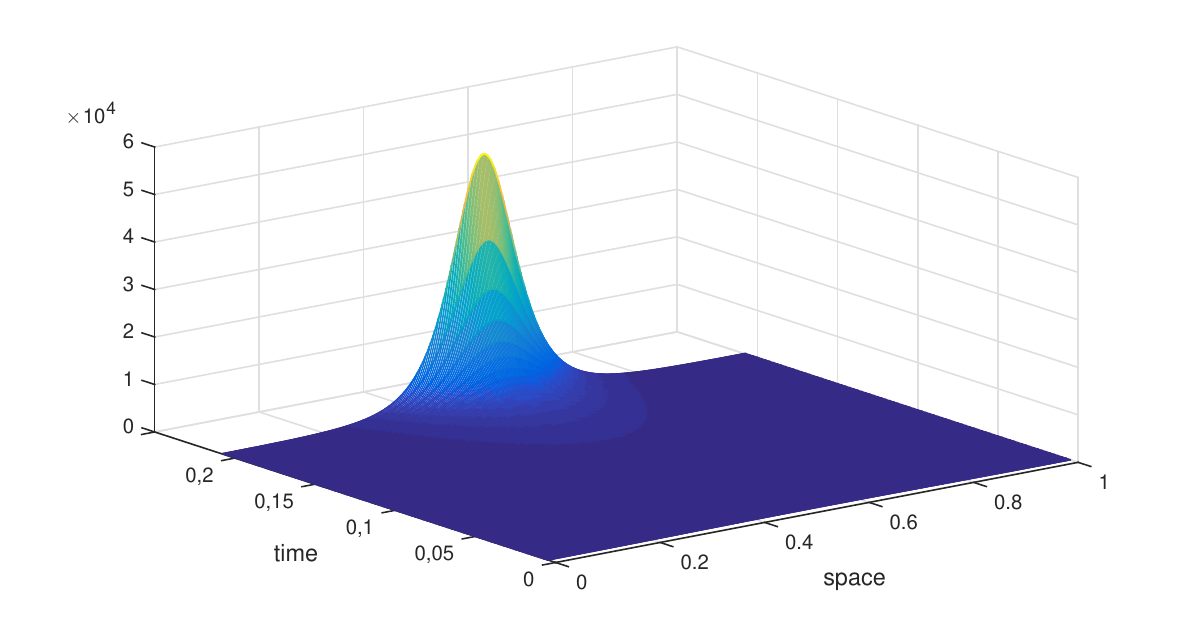}\hspace*{-2mm}
\includegraphics[width=8cm,height=7cm]{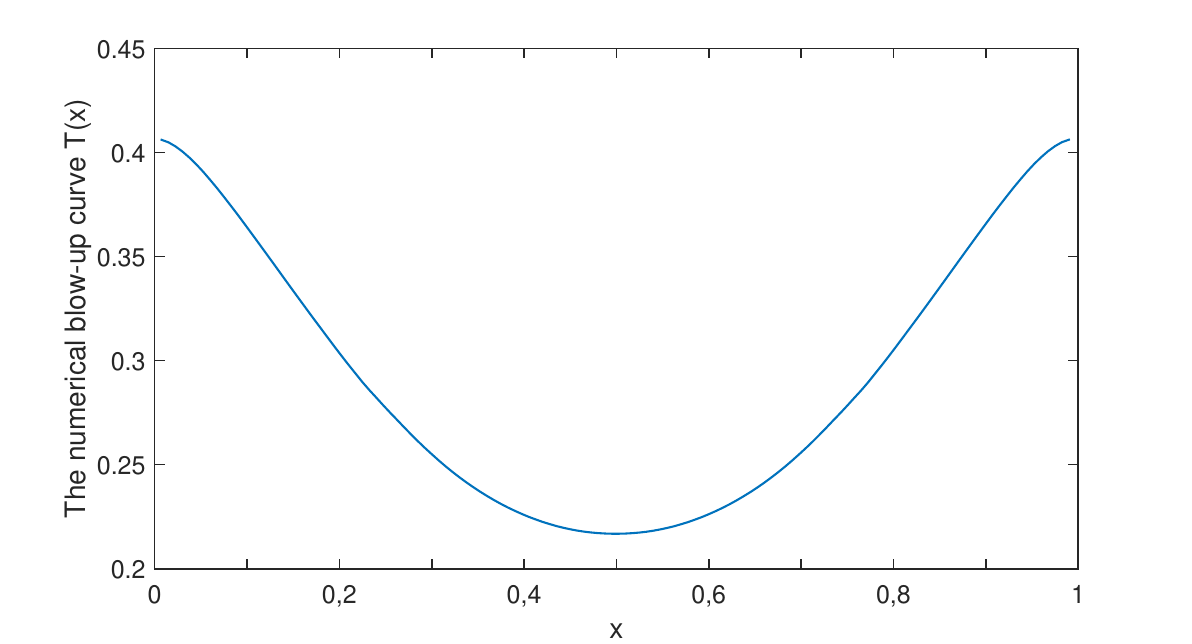}
\caption{\label{numerical solution ex2 power}Left: The numerical solution of example 2. Right: $x$ vs $T(x)$}
\end{center}
\end{figure}

 We consider the system \eqref{equation1} with $p=2$ and the initial data $u_0(x) = 100(1-\cos(2\pi x))$ and $u_1(x) = 10\sin(2\pi x)$. We investigate the numerical blow-up curve by computation of $T^j$ for all $1\leq j\leq J$. Figure \ref{numerical solution ex2 power} shows the numerical solution and the numerical blow up curve. It is shown in \cite{Berg} and \cite{Ngu17} that the value $\tau_k^*$ is independent of $k$ and tends to a constant as $k$ tends to infinity for nonlinear heat equation. We prove that this assertion also holds true in our case. Notice that by \eqref{the solution u(k) power}
\begin{equation*}
u^{(k)}(\xi_k,\tau_k^{*}) = \lambda^{\frac{2}{p-1}}u^{(k-1)}(\lambda\xi_k,\tau_{k-1}^{*}+\lambda\tau^*_k) =\cdots = \lambda^{\frac{2k}{p-1}}u(\lambda^k\xi_k,t_k),
\end{equation*}
where $t_k = \tau^*_{0}+\lambda\tau^*_{1}+\cdots +\lambda^k\tau^*_k$. We recall that if $T$ denotes the blow up time of $u$, then
\begin{align}\label{profile}
(T-t)^{\frac{2}{p-1}}||u(t)||_{\infty} = \mu\  \text{as}\  t\longrightarrow T,\  \text{with}\  \mu = \Big(2\frac{p+1}{(p-1)^2}\Big)^{\frac{1}{p-1}}.
\end{align}   
In particular, at time $t=t_k$, we have
\begin{align*}
(T-t_k)^{\frac{2}{p-1}}||u(t_k)||_{\infty} &= (T-t_k)^{\frac{2}{p-1}}\lambda^{\frac{-2k}{p-1}}||u^k(\tau_k^*)||_{\infty}\\
&= (T-t_k)^{\frac{2}{p-1}}\lambda^{\frac{-2k}{p-1}} M,
\end{align*} 
yielding
\begin{align*}
T-t_k = \lambda^{k}M^{\frac{1-p}{2}} \mu^{\frac{p-1}{2}} + o(1).
\end{align*}
Then, we obtain 
\begin{align*}
\tau_k^* &= \lambda^{-k}(t_k-t_{k-1})\nonumber\\
 &= \lambda^{-k}((T-t_{k-1})-(T-t_k))\nonumber\\
 &= M^{\frac{1-p}{2}}\mu^{\frac{p-1}{2}}(\lambda^{-1}-1) + o(1).
\end{align*}
Finally,
\begin{align}\label{tau}
\lim_{k\longrightarrow\infty}\tau^*_k = M^{\frac{1-p}{2}}\mu^{\frac{p-1}{2}}(\lambda^{-1}-1).
\end{align}
The values of $\tau^*_k$ are tabulated in Table \ref{tab_tau} for various values of $k$. These experimental results shows that $\tau^*_k$ tends to the constant indicated in \eqref{tau} as $k$ tends to infinity, and are in total agreement with our theoretical study. 
\begin{table}[h!]
\begin{center}
\begin{tabular}{|l|c|r|r|r|}
\hline 
$k$  &   $I=100$ &  $I=200$ &  $I=300$ & $I=400$\\ 
\hline 
$10$  & $0.0840$ & $ 0.0855$ & $0.0858$ & $ 0.0860$\\ 
\hline 
$20$  & $  0.0840$ & $0.0855$ & $0.0858$ & $ 0.0860$\\ 
\hline 
$30$  & $0.0840$ & $0.0855$ & $0.0858$ & $ 0.0860$\\
\hline
$40$ & $0.0841$  & $0.0855$ & $ 0.0859$ & $ 0.0862$\\
\hline      
\end{tabular}
\end{center}
\caption{\label{tab_tau} Various values of $\tau^*_k$ with $p=2$.}
\end{table}
%
%
\subsection*{Example 3}\label{example 3}
\begin{figure}
\begin{center}
\includegraphics[width=9cm,height=7cm]{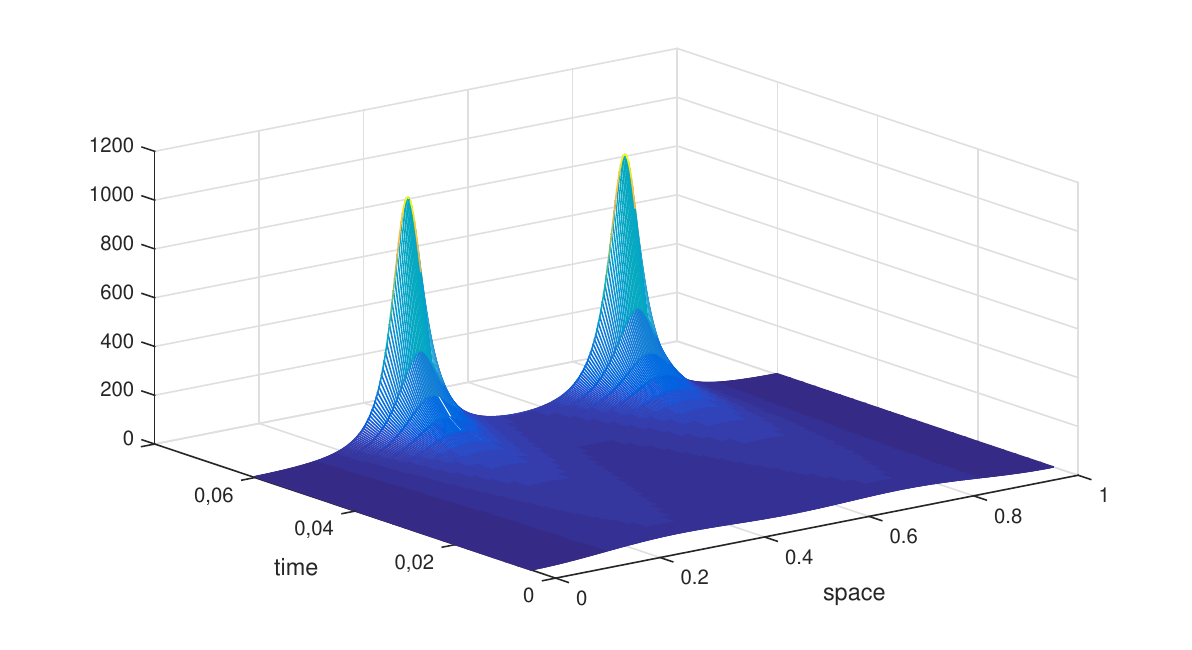}\hspace*{-2mm}
\includegraphics[width=8cm,height=7cm]{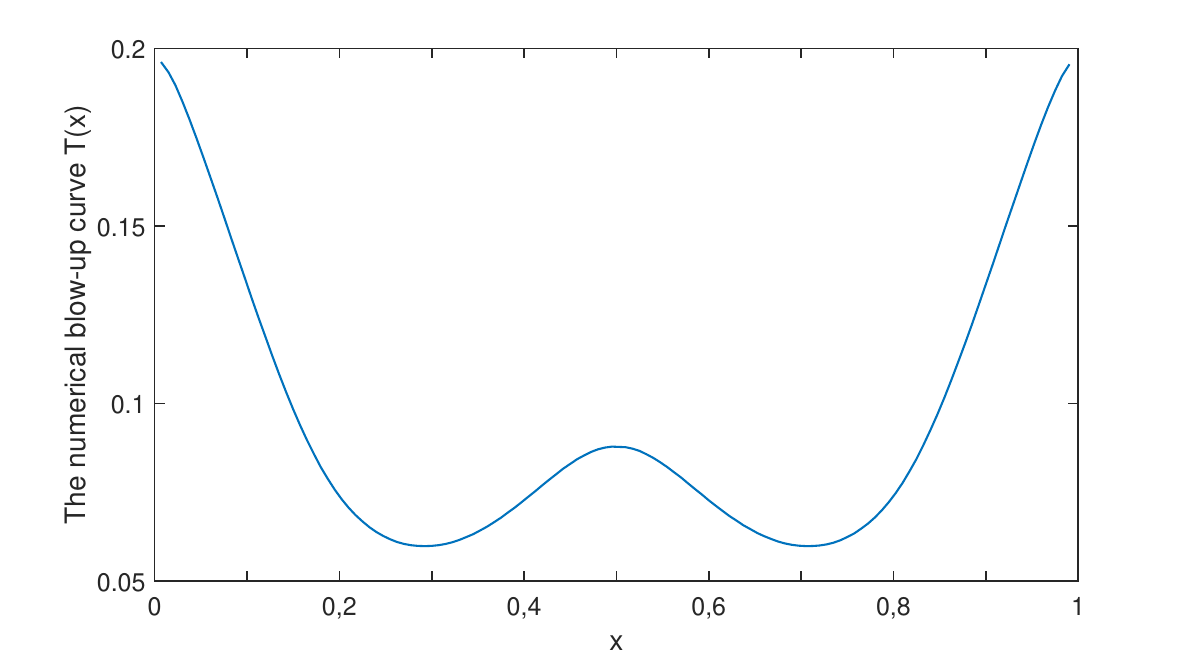}
\caption{\label{numerical solution of exp3}Left: The numerical solution of example 3. Right: $x$ vs $T(x)$.}
\end{center}
\end{figure}

In this example, we consider the system \eqref{equation1} with 
$$u_0(x) = 10(2-\cos(2\pi x)-\cos(4\pi x)), $$
$$u_1(x) = 0,$$
 Figure \ref{numerical solution of exp3} shows the evolution of the numerical solutions in space-time axes for $p = 3$ and the numerical blow-up time $T(x)$.  
\section{Conclusion}
In this paper, we derived a numerical scheme based on both finite difference scheme and a rescaling method for the approximation of the nonlinear wave equation. We proved that under some suitable hypotheses, the numerical solution converges toward the exact solution of the problem. Finally, some numerical experiments are performed and confirm the theoretical study. We expect that all the presented results remain valid for a non linearity $F$  such that $F(u), F'(u),F''(u) \geq 0$ if $u\geq 0$. This will be the object of a future work
\section*{Data Availability Statement}
Data sharing not applicable to this article as no data sets were generated or analyzed during the current study.

\end{document}